\documentclass{elsart1}
\usepackage{mathrsfs}
\usepackage{amsmath}
\usepackage{amssymb}
\input{psfig.sty}

\newtheorem{theorem}{Theorem}[section]
\newtheorem{lemma}{Lemma}[section]

\theoremstyle{definition}

\newtheorem{example}{Example}[section]

\theoremstyle{remark}
\newtheorem{remark}{Remark}[section]

\begin{document}

\begin{frontmatter}

\title{Multiplicity of Periodic Solutions for
Differential Equations Arising in the Study of a Nerve Fiber
Model}

\author[label1,label2]{Chiara Zanini, Fabio Zanolin}
\address[label1]{SISSA-ISAS,
via Beirut 2-4, 34014 Trieste, Italy\\
mailto: zaninic@sissa.it}
\address[label2]{
Dipartimento di Matematica e Informatica, Universit\`a,
via delle Scienze 206,
33100 Udine,
Italy
\\
mailto: zanolin@dimi.uniud.it}

\begin{abstract}
\noindent We deal with the periodic boundary value problem for a
second-order nonlinear ODE which includes the case of the Nagumo
type equation $v_{xx} - g v + n(x) F(v) = 0,$ previously
considered by Grindrod and Sleeman and by Chen and Bell in the
study of the model of a nerve fiber with excitable spines. In a
recent work we proved a result of nonexistence of nontrivial
solutions as well as a result of existence of two positive
solutions, the different situations depending by a threshold
parameter related to the integral of the weight function $n(x).$
Here we show that the number of positive periodic solutions may be
very large for some special choices of a (large) weight $n.$ We
also obtain the existence of subharmonic solutions of any order.
The proofs are based on the Poincar\'{e} - Bikhoff fixed point
theorem.
\end{abstract}


\date{}
\end{frontmatter}

\section{Introduction} \label{zz-sect-1}
This paper deals with the study of periodic solutions to a class of
second order ordinary differential equations arising in the search of stationary
solutions of a partial differential system modelling a nerve
fiber with excitable spines (see \cite{ChBe-94}, \cite{GrSl-85}).

The model equation considered in \cite[pp.391--395]{ChBe-94} takes
the form
\begin{equation}\label{eq-1.1}
v_{xx} - g v + n(x) F(v) = 0
\end{equation}
where $g > 0$ is a given constant, $n(\cdot)$ is a positive
$\beta$-periodic piecewise constant function and $F: {\mathbb R}\to {\mathbb R}$
is a smooth mapping with $-F$ having a $N$-shape.
Typically, $F(s)$ has three zeros $0 < a < 1$ and $F(s) > 0$
for $s\in (-\infty,0[\,\cup\,]a,1[\,,$
$F(s) < 0$
for $s\in \,]0,a[\,\cup\,]1,+\infty).$
According to \cite[pp.391--392]{ChBe-94},
a possible example for $F$ is given by
$F(s) = f({\tilde R}(s)),$ where
$f(u) = u(1-u)(u-a)$ and ${\tilde R}$ is a smooth monotone increasing function
with ${\tilde R}(0) = 0,$ ${\tilde R}(a) = a$ and ${\tilde R}(1) = 1.$
Actually, equation $(\ref{eq-1.1})$ comes from the system
\[
\left\{
\begin{array}{ll}
v_{xx} - g v + n(x) f(v) = 0\\
v = S(u;\rho):= u - \rho^{-1} f(u)\\
\end{array}
\right.
\]
by assuming the function $S(u)$ invertible with inverse ${\tilde R}(v).$
Such an assumption is satisfied when the spine stem resistance $R_S = \rho^{-1}$
approaches zero (cf. \cite[pp.386--387]{ChBe-94} ).

Equations of the form of $(\ref{eq-1.1})$ present a typical threshold phenomenon that
may be easily described if we think for a moment at the weight
$n(x)$ as a constant function and look for constant solutions. It
is clear that if $n = n_0$ is a small positive number, then the
only $\beta$-periodic solution is the trivial one. On the other
hand, if $n = n_1$ is sufficiently large, the line $y = g s$
intersects the {\small$\backslash\!/\!\backslash$}-shaped curve $y = n_1 F(s)$ in two nontrivial
points (at least) and therefore we have the existence of at least
two positive $\beta$-periodic solutions (see Figure 1).

\bigskip

\begin{figure}[ht]
\quad\psfig{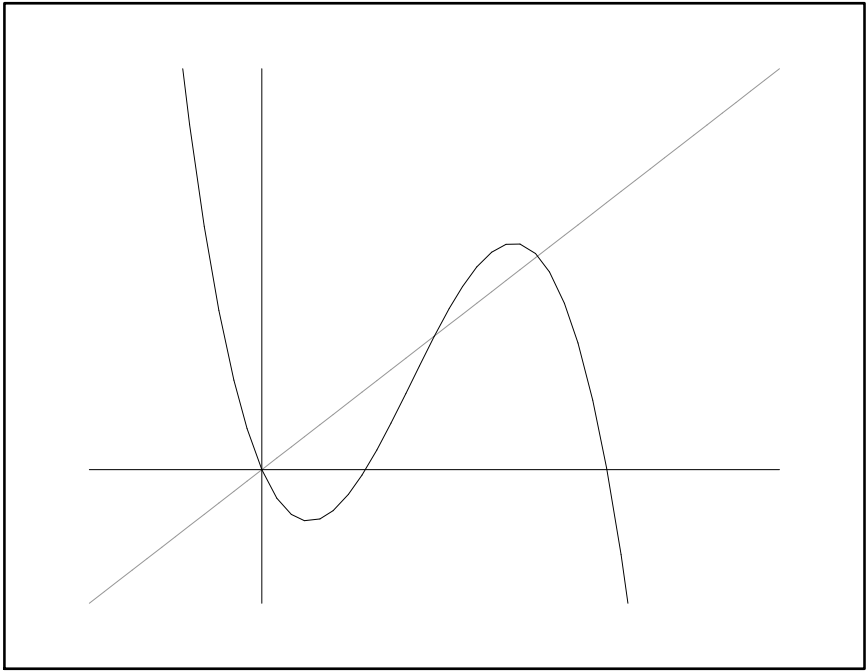}
\vspace*{0in} \caption{\em Intersections of $y=n F(s)$ with $y= gs.$}
\end{figure}

\bigskip

In \cite{ChBe-94} the
authors, using a phase-plane analysis and gluing together the
solutions found for different values of the coefficients,
showed that a similar result can be obtained for a
piecewise constant weight $n(x)$ which, in a period, takes two
values $n_1 > n_0 > 0$ (with $n(x) = n_0$ on $\,]\alpha,\beta[\,$
and $n(x) = n_1$ on $\,]0,\alpha[\,$\,). Then, they obtained
\cite[Lemma 4.1]{ChBe-94} the existence of the only trivial
solution if $\alpha$ and $n_0$ are small (which means that the
weight is close to a small value $n_0$ for the main part of the
time)
and the existence of one positive $\beta$-periodic solution if $n_1$ is
sufficiently large and $\alpha$ is not too close to zero
(which means that the weight is large for an adequate amount of time).
In a further result of the same paper,
the authors \cite[Lemma 4.3]{ChBe-94} also claimed
the existence of ``many'' periodic solutions when $S(u)$ is not invertible
and $n_1$ is large enough. The proof, however, in this
case is only suggested ``by reading the superimposed phase portraits''.
For all these results, some further technical conditions (that we do not
recall here) were assumed. Among such conditions, an hypothesis which implies
\begin{equation}\label{int-1}
\int_0^1 F(s)\,ds > 0,
\end{equation}
was required.

A few questions may arise from these conclusions in
\cite{ChBe-94}, in particular, concerning the possibility of
extending Chen and Bell results to a general class of weight
functions and to the existence of at least two positive
$\beta$-periodic solutions when the weight is sufficiently large
(in some sense), as well as to prove in a more general way (i.e.,
using some argument which is not based by reading the superimposed
phase portraits and therefore applicable only to weights which are
two-step functions or small perturbations of two-steps functions
\cite[Remark 4.6]{ChBe-94}) the existence of ``many'' periodic
solutions for large but arbitrarily shaped profiles.

In a previous recent article \cite{ZaZa-05}, we addressed our attention to the first two questions,
as well as to some related problems. In particular, we showed that the same
features are preserved for a broad class of functions $F$ and for an
arbitrary positive $\beta$-periodic weight $n\in L^1([0,\beta]).$
In this case, the $L^1$-norm (or the mean value) of $n(x)$ plays the role
of a threshold parameter in the sense that we have only the trivial solution
when $|n|_1$ is small and at least two positive $\beta$-periodic solutions
when $|n|_1$ is sufficiently large. Indeed, these results are proved
in \cite{ZaZa-05}  for rather
general second order equations including $(\ref{eq-1.1})$ as a particular case.
The mathematical tools employed in \cite{ZaZa-05} were based on upper and lower solutions, the study of quadratic forms in
Hilbert spaces and critical point theory. Also in \cite{ZaZa-05} condition $(\ref{int-1})$ was
needed for the proof of the existence of two nontrivial solutions.

\smallskip
It is the aim of this work to consider now the third question, that is, to
prove the existence of a large number of periodic solutions
(both harmonics and sub-harmonics), for a general (i.e., not necessarily piecewise constant)
positive
periodic weight function. Note that now the elementary observation made above
which suggested to look at the intersections of the graph of $y= nF(s)$ with the
line $y=gs$ in order to prove the existence of at least two periodic solutions
also for a non-constant weight does not work anymore. Here we have better
to look at the phase-plane portrait of the planar system
\begin{equation}\label{sys-1.1b}
\left\{
\begin{array}{ll}
v' = y\\
y' = g\, v - n(x) F(v)\\
\end{array}
\right.
\end{equation}
and observe that if $n(x) = n$ is a constant function,
with $n$ sufficiently large and $F'(a) > 0$ then, near the point $(a,0)$
an equilibrium point $(a_n,0)$ appears and such an equilibrium point is a local center
surrounded by a family of periodic orbits contained in the strip $]0,1[\,\times
{\mathbb R}.$ The fundamental period of these orbits becomes larger as their energy increases,
 a situation which is reminiscent to the one encountered in the study of the trajectories
 approaching the separatrices in the  nonlinear simple pendulum equation
 (see Figure 2).

\bigskip

\begin{figure}[ht]
\quad\psfig{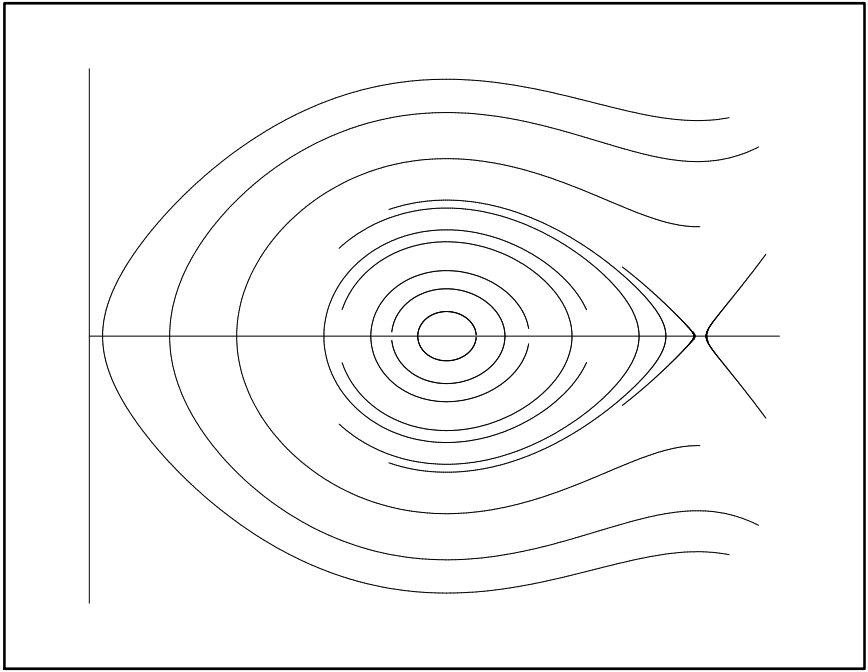}
\vspace*{0in} \caption{\em Phase-portrait of system $(\ref{sys-1.1b})$ for $n(x) = n = 20$
(constant), $g=0.1$ and $F(x) = x(1-x)(x-a),$ with $a=0.6.$ Plotting the orbit-segments at different initial points
$(x_0,0)$ with $x_0 \in \,]0,1],$ but along a same time interval $[-T,T],$ one can see the
presence of periodic trajectories around the stable equilibrium point $(a_n,0)$ which is
near to $(a,0).$ The period needed to complete one turn
becomes larger as the trajectories approach the separatrix.}
\end{figure}

\bigskip

This point of view addresses our investigation toward the search of
periodic trajectories ``near'' $(a,0)$ also when $n(x)$ is not constant,
provided that $n(x)$ has a large principal part.
Indeed, we will show that such periodic trajectories actually exist
(see Theorem \ref{th-5.1} of Section \ref{zz-sect-3}).
More precisely,
using the Poincar\'{e} - Birkhoff fixed point theorem
we shall prove that, if we split
\begin{equation}\label{sn}
n(x) = {\bar n} + {\tilde n}(x),
\end{equation}
where ${\bar n}$ a suitably chosen constant value
(for instance, in some situations, we could take
$\displaystyle{{\bar n}= \frac 1 {\beta} \int_0^{\beta} n(x)\,dx},$
but other possible choices may be suitable as well, depending on the weight function)
then, the number of positive $\beta$-periodic solutions becomes large if
${\bar n}$ grows to infinity and $|{\tilde n}|_1$ is small (in a suitable sense).
Moreover, the same fact is true
also with respect to the subharmonic solutions (that is the $m\beta$-periodic solutions
having $m\beta$ as their minimal period).
To this aim, we shall consider Eq. $(\ref{eq-1.1})$
as a modification of the autonomous equation
\[
v_{xx} - g v + {\bar n} F(v) = 0
\]
and treat the remaining term  ${\tilde n}(x)F(v)$ as a perturbation
(see, for instance, \cite{BuFo-97}, \cite{Di-83}, \cite{FoZa-97}, \cite{HaMa-91} for a similar approach
in the study of some different equations).

As a final remark, we observe that the only crucial assumption for our main multiplicity result of
Section \ref{zz-sect-3} is the existence of a point $a\in \,]0,1[\,$ where
\[
F(a) = 0,\quad F'(a) > 0,
\]
an hypothesis which is always satisfied for the typical $N$-functions considered in the literature.
On the other hand, with our approach, we do not need any condition like $(\ref{int-1})$
on $\int_0^1 F(s)\,ds.$
As pointed out before, we recall that the positivity of the integral in
$(\ref{int-1})$ was required both in the phase-plane analysis approach of Chen and Bell
\cite{ChBe-94} (based on the superposition of the phase portraits for $n = n_0$ and $n = n_1$)
as well as in the variational approach of our recent work
(based on Ambrosetti - Rabinowitz  mountain pass theorem \cite{AmRa-73}).

Hopefully, the study of equation $(\ref{eq-1.1})$ or its variants may be relevant for the
attempt of better understanding the transmission of the impulses along the nerve fibers
also with respect to the fact that there are some serious pathologies
which may arise as a consequence of a deficient distribution of myelin
along the nerve axon (disorders of myelination).
The weight function $n(x)$ in the Nagumo equation represents
the axial distribution of the myelin along the axoplasm.
It is known that some ``fat'' areas alternate with some gaps (the so-called
nodes of Ranvier). Such a profile for the distribution of $n(x)$ has lead
some authors (like in \cite{ChBe-94}, \cite{GrSl-85})
to the study of a piecewise constant weight which alternates
between a small and a large value. However,
the real images of a myelinated axon
show that the true shape of the profile sometimes may be
quite far from this idealized picture and therefore results like
those in \cite{ChBe-94} could not be applied.
Our theorems show that those conclusions are still valid for general positive weights
and {\small$\backslash\!/\!\backslash$}-shaped functions.

\section{Preliminary results and notation}\label{zz-sect-2}
We recall in this section an auxiliary result (taken from {\cite{ZaZa-05}}), about
the periodic solutions of $(\ref{eq-1.1}),$ which is useful for what
follows.

The nerve fiber equation $(\ref{eq-1.1}),$ has the form of a
second-order equation
\begin{equation}\label{eq-1.2}
v'' + h(x,v) = 0
\end{equation}
where $h: {\mathbb R}\times{\mathbb R}\to {\mathbb R}$ is a Carath\'{e}odory function which is
$\beta$-periodic in its first variable, that is,
$h(\cdot,s)$ is measurable for all $s\in {\mathbb R}$
and there is $\beta > 0$ such that
$h(x + \beta,s) = h(x,s)$ for almost every $x\in {\mathbb R}$ and all $s\in {\mathbb R},$
$h(x,\cdot)$ is continuous for almost every $x\in [0,\beta]$ and,
for every $r > 0$ there is a measurable function $\rho_r\in L^1([0,\beta],{\mathbb R}^+)$
such that $|h(x,s)|\leq \rho_r(x)$ for almost every $x\in [0,\beta]$ and
every $s\in [-r,r].$ Solutions of $(\ref{eq-1.2})$ are considered in the Carath\'{e}odory
sense too.

We look for the existence of $m\beta$-periodic solutions to Eq. $(\ref{eq-1.2}),$  for some positive integer
$m,$
that is
for solutions of problem
\[
v'' + h(x,v) = 0,\quad
v(x + m\beta) = v(x),\;\forall\, x\in {\mathbb R},
\leqno{(P_m)}
\]
or, equivalently,
$\displaystyle{
v'' + h(x,v) = 0,}$ with
$v(m\beta) - v(0) = v'(m\beta) - v'(0) = 0.$

Clearly, in the model of our interest, we have
\begin{equation}\label{eq-form}
h(t,s):= - g s + n(t)F(s)
\end{equation}
and the $N$-shape of $-F(s)$ implies that for $h(t,s)$ the
following conditions are satisfied.
\[
h(x,0) \equiv 0
\leqno{(A_1)}
\]
and, for a.e. $x\in [0,\beta],$
\[
h(x,s) > 0, \;\forall \, s < 0,\qquad
h(x,s) < 0, \;\forall \, s \geq 1.
\leqno{(A_2)}
\]
Under these assumptions, we obtain:

\begin{lemma}\label{lem-1}\cite{ZaZa-05} Suppose $(A_1)$ and $(A_2)$ hold.
Then, any possible solution $v(\cdot)$ of problem $(P_m)$ satisfies
$0\leq v(x) \leq 1,\,\forall\, x\in{\mathbb R}.$ Moreover, if $h$ is
(locally) lipschitzean at $s=0,$ then any nontrivial solution of
$(P_m)$ is strictly positive and, if $h$ is (locally) lipschitzean at $s=1,$ then
any solution of $(P_m)$ is strictly less than one.
\end{lemma}

\noindent
The proof is based on some direct estimates and classical arguments from the theory of
upper and lower solutions (see, e.g., \cite{De-95}, \cite{DeHa-96}). We refer to \cite{ZaZa-05}
for more details and remarks.

Thanks to Lemma \ref{lem-1} we can confine ourselves to the solutions $v(x)$
belonging to the interval $[0,1]$
and, moreover, we can modify as we like the function $h(t,s)$ (respectively the
function $F(s)$ in $(\ref{eq-1.1})$ ) on $s\in {\mathbb R}\setminus \,[0,1]\,$
and such a modification will have no effect on the
existence of the periodic solutions as long as the sign conditions outside the
interval $[0,1]$ are preserved.

In the sequel, standard notation is used. We only warn that
for a $\beta$-periodic locally integrable function $u(\cdot),$ we denote by
$|u|_1$ its $L^1$-norm on a interval of length $\beta.$ Sometimes
(when no confusion may occur) the same symbol will be used to denote the $L^1$-norm
on $[0,m\beta],$ when we are interested in the search of $m\beta$-periodic functions.

\section{Multiplicity results}\label{zz-sect-3}
Now we focus our attention on equation $(\ref{eq-1.1})$ and assume that
$g > 0$ is a fixed constant and $F: {\mathbb R}\to {\mathbb R}$ is a continuously differentiable function such that
$
F(0) = F(a) = F(1) = 0
$
with $a: \, 0 < a < 1$ and such that $F(s) > 0$
for $s\in (-\infty,0[\,\cup\,]a,1[\,,$
$F(s) < 0$
for $s\in \,]0,a[\,\cup\,]1,+\infty).$
We also suppose that $n(\cdot): {\mathbb R}\to {\mathbb R}^+$ is a $\beta$-periodic locally integrable function.

We are going to show an application of the Poincar\'{e} - Birkhoff fixed point theorem to the
search of the periodic solutions of $(\ref{eq-1.1}).$ To this aim and for technical reasons it will
be convenient to modify suitably the function $F$ outside the interval $[0,1].$
Precisely, let us set
\[
\delta(s) : = \max\{0,\min\{s,1\}\}
\]
and define
\[
F_0(s) := F(\delta(s)) + k_0 \ell(s),
\]
where $\ell(s)$ is given by
\begin{equation}\label{ell}
\ell(s) = \,
\left\{
\begin{array}{ll}
\exp(1/s),\; &s < 0\\
0,\; & 0\leq s\leq 1\\
-\exp(1/(1-s)),\; & s >1,\\
\end{array}
\right.
\end{equation}
and
\[
0 < k_0 \leq c_0\,:=\max_{s\in[0,1]} |F(s)|.
\]
With these positions, we have that
$F_0 :{\mathbb R}\to {\mathbb R}$ is a bounded (globally) Lipschitz continuous function which is smooth in $\,]0,1[\,,$
$F_0(s) = F(s)$ for all $s\in [0,1],$ $F_0(s) > 0$ for $s < 0$ and $F_0(s) < 0$ for $s > 1.$
Moreover,
\[
\sup_{s\in{\mathbb R}} |F_0(s)| = \max_{s\in{\mathbb R}} |F_0(s)| = \max_{s\in [0,1]} |F(s)| = c_0\,.
\]
A straightforward application of Lemma \ref{lem-1} ensures that for each positive integer $m,$
all the possible nontrivial $m\beta$-periodic solutions of equation
\begin{equation}\label{eq-5.1}
v'' - g v + n(x) F_0(v) = 0
\end{equation}
have range into the open interval $\,]0,1[\,$ and so they are,
indeed, solutions of equation $(\ref{eq-1.1}).$

The fact that after our modification the nonlinearity is now bounded, suggests the
possibility of splitting $n(x)$ as in $(\ref{sn})$ and
writing equation $(\ref{eq-5.1})$ as
\begin{equation}\label{eq-5.2}
v'' - g v + {\bar n} F_0(v) = p(x,v)
\end{equation}
which looks like a perturbation of the autonomous equation
\begin{equation}\label{eq-5.3}
v'' - g v + {\bar n} F_0(v) = 0
\end{equation}
by an external source
\[
p(x,v):= -{\tilde n}(x) F_0(v).
\]
Now we are in position to prove our main multiplicity result.
\begin{theorem}\label{th-5.1}
Assume that
\[
F'(a) > 0.
\]
Then, for every integer $N \geq 1,$ there exists a (large) value $\mu = \mu_N > 0$ such that
for each $\mu_2 \geq \mu_1 > \mu$ there is a (small) value
$\varepsilon = \varepsilon_{N,\mu_1,\mu_2} > 0$ such that for each $n(\cdot)$ with
\[
{\bar n} \in [\mu_1,\mu_2]\quad\mbox{and }\; |{\tilde n}|_1  < \varepsilon
\]
there are at least $2N$ solutions of equation $(\ref{eq-1.1})$ which are $\beta$-periodic
and take values in $\,]0,1[\,.$
\end{theorem}

The proof of Theorem \ref{th-5.1} is based on the Poincar\'{e} -  Birkhoff fixed
point theorem and will be split into some steps.

First of all, we write equation $(\ref{eq-5.1})$ as a first order planar system
of the form
\begin{equation}\label{sys-5.1}
\left\{
\begin{array}{ll}
\dot x = y\\
\dot y = - h(t,x)\\
\end{array}
\right.
\end{equation}
with
\[
h(t,s) = - g s + n(t)F_0(s),\quad \dot x = \frac{dx}{dt}.
\]
Note that we have changed the name to the
variables, by denoting the $x$ and the $v$ variables in $(\ref{eq-5.1})$ as
$t$ and $x,$ respectively. According to our hypotheses, $h: {\mathbb R}\times{\mathbb R}\to {\mathbb R}$
is a Carath\'{e}odory function which is $\beta$-periodic in $t$ and satisfies
a global Lipschitz condition
\[
|h(t,s_1) - h(t,s_2)| \leq A(t) |s_1 - s_2|, \quad \mbox{for all } s_1\,, s_2\, \in {\mathbb R} \;
\mbox{and for a.e. } t\in {\mathbb R},
\]
where $A(\cdot)$ is a suitable $\beta$-periodic measurable function with $A(\cdot) \in L^1([0,\beta]).$
In our case, we actually have $A(t) = g + n(t)L_0$ where $L_0$ is a Lipschitz constant for $F_0$
(which is globally lipschitzean on ${\mathbb R}$). As a
consequence, we have that for each $z_0 = (x_0,y_0)\in {\mathbb R}^2$ and
$t_0\in {\mathbb R},$ system $(\ref{sys-5.1})$ has a unique solution
$\zeta(t) = \zeta(t,t_0,z_0) = (x(t,t_0,z_0),y(t,t_0,z_0))$
satisfying the initial condition $\zeta(t_0) = z_0\,,$ with
$\zeta(t)$ defined for all $t\in {\mathbb R}.$ Hence, the Poincar\'{e}'s
map
\[
\phi: z_0 \mapsto \zeta(\beta,0,z_0)
\]
is well defined as a homeomorphism of ${\mathbb R}^2$ onto itself.

Next, we introduce a rotation number associated to $\phi$ with respect to a given point
$q_0\in{\mathbb R}^2.$

Let $m\in {\mathbb N}$ be a positive integer and let
$q_0 = (q^0_1,q^0_2)$ and  $z_0 = (x_0,y_0)$ be such that
\[
\zeta(t,0,z_0)\not= q_0\,\quad\forall\, t\in [0,m\beta].
\]
We define the rotation number (see, for instance, \cite{DiZa-93}, \cite{Za-03}) as
\[
\mbox{rot}_m(z_0,q_0) = \frac{1}{2\pi}\int_0^{m\beta}
{
\frac{y(t)^2 + x(t)h(t,x(t))}{(x(t) - q^0_1)^2 + (y(t)- q^0_2)^2}\,dt
}
\]
which represents the normalized angular displacement around the point $q_0$
of the solution $\zeta(t,0,z_0)$ for $t$ varying along the time interval $[0,m\beta].$
In fact, if we use the Pr\"{u}fer transformation and write the solution
$\zeta(t)$ in polar coordinates $(\theta,\rho)$ with respect to the point $q_0\,,$
we have that $\rho(t)^2 = ||\zeta(t,0,z_0) - q_0||^2 = (x(t,0,z_0) - q^0_1)^2 + (y(t,0,z_0)- q^0_2)^2 > 0$
for all $t\in [0,m\beta]$ and therefore the number $\theta(t) - \theta(0)$ is well defined.
It turns out that
\[
\mbox{rot}_m(z_0,q_0) = \frac{\theta(0) - \theta(m\beta)}{2\pi}
\]
and this number counts as positive the clockwise rotations around the point $q_0\,.$

Our main tool to prove the existence of periodic solutions is the
following result which is adapted to our situation
from W. Ding's generalized version of
the Poincar\'{e} -  Birkhoff theorem \cite{Di-83}.

\begin{theorem}\label{th-5.2}
Let $q_0=(q^0_1,q^0_2)\in{\mathbb R}^2$ and let $D_0$ be an open neighborhood of $q_0$ such that
\begin{equation}\label{eq-nonq}
\zeta(t,t_0,z_0) \not=q_0\,,\quad\forall\, z_0\in \partial D_0\,,\;\;\forall\, t_0\in[0,m\beta[\,,\;\forall\, t\in [t_0,m\beta].
\end{equation}
Suppose that $\Gamma\subset {\mathbb R}^2\setminus D_0$ is a simple closed curve which is star-shaped with respect to $q_0$
and there exist a positive integer $j$
such that
\begin{equation}\label{eq-rot1}
\mbox{\rm rot}_m(z_1,q_0) > j,\quad\forall\, z_1\in \Gamma.
\end{equation}
Furthermore, let us assume there is a (sufficiently large) radius
$R > 0,$ with $\Gamma \subseteq B(q_0,R)$ such that
\begin{equation}\label{eq-rot2}
\mbox{\rm rot}_m(z_2,q_0) < 1,\quad\forall\, ||z_2|| = R.
\end{equation}
Define ${\mathcal A}$ to be the open annulus bounded by $\Gamma$ and $\partial B(q_0,R).$
Then, there are at least $2j$ solutions ${\tilde{\zeta}}_k$ and
${\hat{\zeta}}_k$ (for $k=1,\dots,j$) of $(\ref{sys-5.1})$ which are $m\beta$-periodic
and such that ${\tilde{\zeta}}_k(0), {\hat{\zeta}}_k(0) \in {\mathcal A}$ and
\[
\mbox{\rm rot}_m(z_0,q_0) = k,\quad\mbox{for }\; z_0 = {\tilde{\zeta}}_k(0)\;
\mbox{or }\,
z_0 = {\hat{\zeta}}_k(0).
\]
In the particular case when $q_0 = (q^0_1,0)$ we have also that, setting
${\tilde{\zeta}}_k(t) = ({\tilde{u}}_k(t),{\tilde{y}}_k(t))$ and
${\hat{\zeta}}_k(t) = ({\hat{u}}_k(t),{\hat{y}}_k(t)),$
it follows that ${\tilde{u}}_k$ and ${\hat{u}}_k$ are $m\beta$-periodic solutions
of equation
\[
u'' + h(t,u) = 0
\]
with ${\tilde{u}}_k(t) - q^0_1$ and
${\hat{u}}_k(t) - q^0_1$ having precisely $2k$ simple zeros in the interval $[0,m\beta[\,.$
\end{theorem}

\noindent
{\em Proof. }
Let $\phi :{\mathbb R}^2\to{\mathbb R}^2$ be the Poincar\'{e}'s map associated to system $(\ref{sys-5.1}).$
By the Liouville theorem it follows that $\phi$ is area-preserving and therefore
$\phi^m$ is an area-preserving planar homeomorphism as well.
Let $E_0$ be the open bounded set confined by the curve $\Gamma.$ By the assumption,
$q_0\in D_0\subset E_0$ and $\partial E_0 = \Gamma.$
Note also that the rotation number $\mbox{\rm rot}_m(z_0,q_0)$ is well defined for
each $z_0\in {\mathbb R}^2\setminus E_0$ due to the fact $\zeta(t,0,z_0)\not= q_0$ for all $t\in [0,m\beta]$
and each $z_0\not\in E_0\,.$
Indeed, if, by contradiction, there is $z_0\not\in E_0$ such that
$\zeta(t^*) = \zeta(t^*,0,z_0) = q_0$ for some $t^*\in \,]0,m\beta],$ then $\zeta(t_1,0,z_0) = z_1\in\partial D_0$
for some $t_1\in [0,t^*[\,.$ This means that $\zeta(t^*) = \zeta(t^*,t_1,z_1) = q_0$ and
we have a contradiction to hypothesis $(\ref{eq-nonq}).$
Now, if we restrict $\psi := \phi^m$ to the closed disk $B[q_0,R]$ we have that
\[
\psi^{-1}(q_0) \in E_0
\]
and $\psi$ satisfies the twist condition with respect to the inner and the outer boundaries of
the closed annulus ${\bar{\mathcal A}}$ for each integer $k\in [1,j].$ Namely, we have
$\mbox{\rm rot}_m(\cdot,q_0) > k,$ on $\Gamma$ and
$\mbox{\rm rot}_m(\cdot,q_0) < k,$ on $\partial B(q_0,R).$
The
W. Ding's generalized version of the twist theorem guarantees the existence of at least
two fixed points in ${\mathcal A}$ for $\psi.$ These fixed points are initial values
(at the time $t=0$) of two $m\beta$-periodic solutions ${\tilde{\zeta}_k}(\cdot)$
and ${\hat{\zeta}}_k(\cdot)$ of system $(\ref{sys-5.1}),$ respectively.
Moreover, the rotation number associated to ${\tilde{\zeta}_k}$ and ${\hat{\zeta}}_k$
is equal to $k.$

In the special case when $q_0 = (q^0_1,0),$ and
by virtue of the particular form of system $(\ref{sys-5.1}),$
we know that the result about the rotation numbers implies
that the first coordinates ${\tilde{u}_k}$ and ${\hat{u}}_k$ of ${\tilde{\zeta}_k}$ and ${\hat{\zeta}}_k$
crosses the line $u = q^0_1$ exactly $2k$ times in the interval $[0,m\beta[\,.$
\hfill$\square$

\begin{remark}\label{rem-5.1}
{\rm
Clearly, if $u(\cdot)$ is a $m\beta$-periodic solution of the equation $u'' + h(t,u) =0$
with $h(t+\beta,u) = h(t,u)$ then
\[
u(\cdot + \beta), \dots, u(\cdot + j\beta),
\dots, u(\cdot + (m-1)\beta)
\]
are $m\beta$-periodic solutions as well. We consider these solutions as equivalent
each other
and we say that they belong to the same periodicity class.
A further consequence of the Poincar\'{e} - Birkhoff theorem (as pointed out in \cite{Ne-77})
ensures that the solutions ${\tilde{u}_k}$ and ${\hat{u}_k}$
(for which the existence is claimed in Theorem \ref{th-5.2}
do not belong to the
same periodicity class. Obviously, also the solutions
${\tilde{u}_k}\,,$ ${\hat{u}_k}$
${\tilde{u}_{\ell}}\,,$ ${\hat{u}_{\ell}}$ for $\ell\not=k$
belong to different periodicity classes
(in fact their rotation numbers are different). We refer also to \cite{ReZa-96}, \cite{Za-98}
for a throughout discussion concerning this aspect.
\\
Moreover, we observe that the information about the associated rotation numbers permits obtain
some conclusions about the minimality of the period.
For instance, if $m\geq 2$ and $k\geq 1$ are co-prime numbers (that is $m/k$ is not further reducible),
then it is possible to prove that the $m\beta$-periodic solutions ${\tilde{u}_k}$ and ${\hat{u}_k}$
are not $i\beta$-periodic for each $i=1,\dots,m-1.$ In particular, for $k=1,$ we have that
the solutions we find have $m\beta$ as their minimal period (see \cite{DiZa-93}
where this discussion is carried on with more details).
}
\end{remark}

In the proof of Theorem \ref{th-5.1} using Theorem \ref{th-5.2} we take
\[
q_0 = (a_{\bar{n}},0)
\]
with $0 < a_{\bar{n}} < 1$ such that
\[
\frac{F_0(a_{\bar{n}})}{a_{\bar{n}}} = \frac{g}{\bar n}
\]
Note that
\[
a_{\bar{n}} > a,\quad\forall\, \bar n > 0
\quad
\mbox{and }\;
\lim_{\bar n \to +\infty} a_{\bar{n}} = a.
\]

Our next result shows that large solutions rotate slowly.

\begin{lemma}\label{lem-5.1}
There is $R_0 = R_0(m,|n|_1)> 0$ such that for each initial point $z_2$ with $||z_2|| \geq R_0$
\[
\mbox{\rm rot}_m(z_2,q_0) < 1,\quad\forall\, ||q_0|| \leq 1.
\]
\end{lemma}

\noindent
{\em Proof. }
First of all, we recall a well known consequence of the global existence of the solutions
(cf. \cite{Kr-68}), that is,
{\em there is a continuous nondecreasing function $\eta: {\mathbb R}^+\to {\mathbb R}^+_0\,,$ with $\eta(r) > r$
for all $r\geq 0,$ such that
\[
||\zeta(t,t_0,z_0)|| > R,\quad \forall\, t,t_0 \in [0,m\beta],\; \forall\, ||z_0||\geq \eta(R).
\]
}
The function $\eta$ depends upon the constants which bound the growth of $h(t,x)$
and hence, we have $\eta = \eta_{m,|n|_1}$ (see also \cite[Lemma 2]{Za-96}).

Let $R > 0$ be such that $||\zeta(t)||\geq R$ for all $t\in [0,m\beta].$ We claim that
if $R > 1$  then
\[
\mbox{\rm rot}_m(\zeta(0),0) < 1/2.
\]
Indeed, let us assume the contrary and suppose that $\mbox{\rm rot}_m(\zeta(0),0) \geq 1/2.$
This implies that the projection $\zeta(t)/||\zeta(t)||$ of $\zeta(t)$ to $S^1$ covers at least
one of the intersections of $S^1$ with the standard quadrants in the $x\,y$-plane.
Just to fix the ideas, let us assume that the trajectory crosses the first quadrant.
Hence, there are $0\leq t_1 < t_2 \leq m\beta$ such that
\begin{itemize}
\item[]
$\displaystyle{
\zeta(t_1) = (x(t_1),y(t_1)) = (1,y_1),\quad \mbox{with }\; y_1 \geq (R^2 - 1)^{\frac 1 2}
}$
\item[]
$\displaystyle{
\zeta(t_2) = (x(t_2),y(t_2)) = (x_2,0),\quad \mbox{with }\; x_2 \geq R
}$
\item[]
$\displaystyle{
||\zeta(t)|| \geq R,\quad x(t) > 1,\quad y(t) > 0,\quad\forall\, t\in \, ]t_1,t_2[\,.
}$
\end{itemize}
Then, passing to the polar coordinates (via the Pr\"{u}fer transformation), we find, for
$t\in [t_1,t_2]$ and using the fact that $h(t,s)s \leq 0,$ for all $s\geq 1,$
\[
- \theta'(t) = \frac{y(t)^2 + x(t)h(t,x(t))}{x(t)^2 + y(t)^2}\leq \frac{y(t)^2}{x(t)^2 + y(t)^2}
= \sin^2(\theta(t)).
\]
For $t\in [t_1,t_2[\,,$ we have
\begin{eqnarray*}
\cot(\theta(t)) - \cot(\theta(t_1))
&=&
\int_{\theta(t)}^{\theta(t_1)}\frac{du}{\sin^2 u} =
-\int_{t_1}^t \frac{\theta'(s)}{\sin^2(\theta(s))}\,ds\\
&\leq& \int_{t_1}^t 1\,ds = t-t_1 \leq m\beta.\\
\end{eqnarray*}
Therefore, we get
\[
\cot(\theta(t)) - \frac{1}{(R^2 - 1)^{\frac1 2}} \leq \cot(\theta(t)) - y_1^{-1} \leq m\beta,\quad
\forall\, t\in [t_1,t_2[\,.
\]
This yields to a contradiction as $t\to t_2^-$ and therefore our claim is proved.

\medskip

Adapting a result in \cite[Lemma 2.2]{DiZa-93} to our situation, we
have that for every $\varepsilon > 0$ there is $R^{\varepsilon} > 1$ such that
\[
|\mbox{\rm rot}_m(\zeta(0),0) - \mbox{\rm rot}_m(\zeta(0),q_0)| < \varepsilon,
\]
holds for every point $q_0 = (a_{\bar{n}},0)$ with $||q_0||\leq 1$ and all solutions $\zeta(\cdot)$
with  $||\zeta(t)|| \geq R^{\varepsilon}\,,\forall\, t\in [0,m\beta].$

Thus there exists a sufficiently large radius $R^* > 1$ such that
\[
\mbox{\rm rot}_m(\zeta(0),q_0) < \frac 3 4,\quad
\forall ||q_0||\leq 1,\forall \zeta \;
\mbox{with } ||\zeta(t)|| \geq R^*\,,\forall\, t\in [0,m\beta]
\]
and, by the result recalled at the beginning of the proof, we achieve the
conclusion by taking
\[
R_0 \geq \eta(R^*).
\]
Thus the proof is complete.
\hfill$\square$

As a next step, we evaluate the rotation numbers of small solutions around $q_0 = (a_{\bar n},0).$
To this end, we consider autonomous equation $(\ref{eq-5.3})$ that we write as a first order
system
\begin{equation}\label{sys-5.2}
\left\{
\begin{array}{ll}
\dot x = y\\
\dot y = g x - {\bar n}F_0(x).\\
\end{array}
\right.
\end{equation}
System $(\ref{sys-5.2})$ is a conservative system with energy
\[
E_{\bar n}(x,y) = \frac 1 2 y^2 - \frac 1 2 g x^2 + {\bar n}{\mathcal F}_0(x),\quad
\mbox{with }\; {\mathcal F}_0(x) = \int_0^x F_0(s)\,ds.
\]
For simplicity in the notation, in the sequel we set $E_{\bar n} = E.$
Let us fix $d_0 >0$ and $b\in \,]a,1[\,$ such that
\[
F'_0(x)\geq d_0\,,\quad\forall\, x\in [a,b]
\]
and take
\[
\mu_0 \geq \frac{g}{d_0}
\]
such that
\[
a < a_{\bar n} < b,\quad\mbox{for }\;
{\bar n} > \mu_0\,.
\]
For ${\bar n} > \mu_0$ we have that $E'(\cdot,0)$ is strictly increasing on $[a,b]$ and
vanishes in $a_{\bar n}\,.$
Hence, on the interval $[a,b]$ we have that $E(\cdot,0)$ is a strictly convex function
having absolute minimum at the point $a_{\bar n}\,.$
If we take now a constant $c_{\bar n}$ with
\[
E(a_{\bar n},0) < c_{\bar n} \leq \min\{E(a,0),E(b,0)\},
\]
we have that the level line
\[
\Gamma_{\bar n} = \{(x,y): a\leq x\leq b,\; E(x,y) = c_{\bar n}\}
\]
is a star-shaped curve around the equilibrium point $q_0\,.$
Moreover, for any point $q\in \Gamma_{\bar n}$ we have that the solution of system
$(\ref{sys-5.2})$ with $(x(0),y(0)) = q$ is periodic and the corresponding orbit
coincides with $\Gamma_{\bar n}\,.$

In order to avoid misunderstanding with the previously defined rotation number,
we denote by $\mathscr{R}\!{\it ot}_m$ the rotation number associated to $(\ref{eq-5.3})$
along the time interval $[0, m\beta].$
If we denote by $\tau_{\bar n}$ the fundamental
(i.e., minimal) period of $\Gamma_{\bar n}\,,$ then we can conclude that
\[
\mathscr{R}\!{\it ot}_m(q,q_0) \geq\,\left\lfloor \frac{m \beta}{\tau_{\bar n}}\right\rfloor,
\quad\forall\, q\in \Gamma_{\bar n}\,.
\]

We claim that $\displaystyle{
\lim_{{\bar n}\to + \infty} \tau_{\bar n} = 0.
}$

Indeed, by the well known time-mapping formula we have
\[
\tau_{\bar n} = \tau^-_{\bar n} + \tau^+_{\bar n} =
\sqrt{2}\,\int_{b^-_{\bar n}}^{a_{\bar n}} \frac{d u}{(c_{\bar n} - E(u,0))^{\frac 1 2}}
+
\sqrt{2}\,\int_{a_{\bar n}}^{b^+_{\bar n}} \frac{d u}{(c_{\bar n} - E(u,0))^{\frac 1 2}}
\]
where
\[
b^-_{\bar n} < a_{\bar n} < b^+_{\bar n},\quad E(b^-_{\bar n},0) = E(b^+_{\bar n},0) = c_{\bar n}\,.
\]
For $u\in [a_{\bar n},b^+_{\bar n}],$ we have
\begin{eqnarray*}
c_{\bar n} - E(u,0)
&=& E(b^+_{\bar n},0) - E(u,0) = \int_{u}^{b^+_{\bar n}} E'(x,0)\,dx\\
&=& \int_{u}^{b^+_{\bar n}} (E'(x,0) - E'(a_{\bar n},0))\,dx
= \int_{u}^{b^+_{\bar n}} \left( \int_{a_{\bar n}}^x E''(\xi,0)\,d\xi\right)\,dx\\
&\geq&
\frac{{\bar n} d_0 - g}{2}\,
\left((b^+_{\bar n} - a_{\bar n})^2 - (u - a_{\bar n})^2\right).
\end{eqnarray*}
Therefore, via an elementary integration, we obtain
\[
\tau^+_{\bar n} \leq \frac{\pi}{\sqrt{{\bar n} d_0 - g}}\,.
\]
A similar computation for $\tau^-_{\bar n}$ yields to
\[
\tau_{\bar n}\leq \frac{2\pi}{\sqrt{{\bar n} d_0 - g}}
\]
and this proves the claim.
Hence we have that
\[
\mathscr{R}\!{\it ot}_m(q,q_0) \geq\,\left\lfloor
\frac{m \beta \sqrt{{\bar n} d_0 - g}}{2\pi}\right\rfloor,
\quad\forall\, q\in \Gamma_{\bar n}\,,
\]
which shows that $\mathscr{R}\!{\it ot}_m(q,q_0)\to +\infty$ with the order of $\sqrt{\bar n}.$

\bigskip

Now we are in position to complete the proof of Theorem \ref{th-5.1}.

Fix $m\geq 1$ and $N\geq 1.$
Let $\mu = \mu_{N}\geq \mu_0$ be such that
\[
\mathscr{R}\!{\it ot}_m(q,q_0) = \sigma_{\bar n}\geq \sigma > N, \quad \forall\, q\in \Gamma_{\bar n}\,,
\]
holds for each ${\bar n} > \mu.$
Fix also an interval $[\mu_1,\mu_2]\subset \,]\mu,+\infty)$ and consider equation
$(\ref{eq-5.2})$ with ${\bar n}\in [\mu_1,\mu_2].$ By the continuous dependence of the solutions
from the coefficients and from initial data (otherwise, some direct estimates may also be
performed) we have that there exists $\varepsilon_1 > 0$ (depending on $\mu_1$ and $\mu_2$ and,
in turns, also on $N$) such that for each forcing term $p$ with $|p|_1 < \varepsilon_1\,,$
condition $(\ref{eq-nonq})$ is satisfied as well as
$\mbox{\rm rot}_m(z_0,q_0) > N$ holds for each $q_0 = (a_{\bar n},0)$
and $z_0 \in \Gamma_{\bar n}$ provided that ${\bar n} \in [\mu_1,\mu_2].$
Recalling the definition of $c_0$ as a bound for $|F_0|,$ we have that for
\[
|{\tilde n}|_1 < \varepsilon = \frac{\varepsilon_1}{c_0}\,,
\]
it follows that $(\ref{eq-nonq})$ holds and
\[
\mbox{\rm rot}_m(z_1,q_0) > N, \quad \forall\, z_1\in \Gamma_{\bar n}.
\]
For any chosen $n(x)$ with ${\bar n}\in [\mu_1,\mu_2]$ and $|{\tilde n}|_1 < \varepsilon,$
we can take a sufficiently large radius $R_n$ such that
$\Gamma_{\bar n} \subseteq B(q_0,R_n)$ and
\[
\mbox{\rm rot}_m(z_2,q_0)  < 1, \quad \forall\, ||z_2|| = R_n\,.
\]
The Poincar\'{e} - Birkhoff theorem (Theorem \ref{th-5.2}) ensures the existence of $2N$ solutions
which are $m\beta$-periodic.  Clearly, for $m = 1$ we have exactly the result claimed in
Theorem \ref{th-5.1} and thus the proof is complete.
\hfill$\square$

Clearly, from the above proof we have a result about the existence of subharmonic solutions
(of period $m\beta$ with $\beta \geq 2$) as well. Indeed, the discussion about
the minimality of the period given in Remark \ref{rem-5.1} yields to:
\begin{theorem}\label{th-5.3}
Assume that
\[
F'(a) > 0.
\]
Then, for every integer $m\geq 2$ and each $K \geq 1,$ there is a (large) value $\mu = \mu_K > 0$ such that
for each $\mu_2 \geq \mu_1 > \mu$ there is a (small) value
$\varepsilon = \varepsilon_{K,\mu_1,\mu_2} > 0$ such that for each $n(\cdot)$ with
\[
{\bar n} \in [\mu_1,\mu_2]\quad\mbox{and }\; |{\tilde n}|_1  < \varepsilon
\]
there are at least $2K$ solutions of equation $(\ref{eq-1.1})$ which are $m\beta$-periodic
and take values in $\,]0,1[\,.$ Moreover, such solutions do not belong to the same periodicity
class and also they are not $i\beta$-periodic, for each $i=1,2,\dots, m-1,$ so that
the period $m\beta$ is minimal.
\end{theorem}

\noindent
{\em Proof. }
Once we have fixed $m$ and $K,$ we take the set
\[
{\mathcal N}:=\{l_1=1 < l_2 <\dots < l_K:=N\}
\]
made by the first $K$ numbers
which are co-prime with $m.$\footnote{Two positive integers $l$ and $m$ are co-prime
(or relatively prime) if $\mbox{\rm GCD}(l,m) = 1$}
According to the conclusion of the proof of Theorem \ref{th-5.1},
we find that there are $2N\geq 2K$ periodic solutions of period $m\beta.$ As a consequence of Remark
\ref{rem-5.1} we know that these $2N$ solutions appears in pairs which are characterized by the
fact that they have, respectively, $2, 4, 6, \dots, 2N$ oscillations in the time interval $[0,m\beta[\,.$
Such oscillations are associated to the rotation numbers of the solutions (an information which
comes from the use of the Poincar\'{e} - Birkhoff fixed point theorem).
As we already explained in Remark \ref{rem-5.1} we are sure that for those solutions having
a rotation number $\ell$ which is co-prime with $m,$ it holds that $m\beta$ is the minimal period.
By the choice of $N$ and the set ${\mathcal N}$ we are lead to conclude that at least $2K$
among the $2N$ periodic solutions (that is those with rotation numbers equal to
$l_1\,,l_2\,,\dots, l_K$) are of minimal period. This observation concludes our proof.
\hfill$\square$

\bigskip

We end our paper, by showing that it is possible to realize a ``good'' splitting of the weight function
$n(x)$ for a profile like the one considered in \cite{ChBe-94}.

\begin{example}\label{ex-3.1}
{\rm
Consider now equation $(\ref{eq-1.1})$ and assume, like in
\cite{ChBe-94}
that $n(x)$ is a periodic piecewise constant function satisfying
\begin{equation}\label{enne}
n(x) = \;
\left\{
\begin{array}{ll}
n_1\quad \mbox{if }\, x\in \,]0,\alpha[\, \mod \beta\\
n_0\quad \mbox{if }\, x\in \,]\alpha,\beta[\, \mod \beta\\
\end{array}
\right.
\end{equation}
with $0 < n_0 < n_1\,.$ Like in \cite{ChBe-94} we assume that
$F$ is a smooth function satisfying the sign conditions we have already described at the
beginning, that is
$
F(0) = F(a) = F(1) = 0
$
with $a: \, 0 < a < 1$ and such that $F(s) > 0$
for $s\in (-\infty,0[\,\cup\,]a,1[\,,$
$F(s) < 0$
for $s\in \,]0,a[\,\cup\,]1,+\infty).$ In order to apply our result, we also suppose that
$F'(a) > 0$ (an hypothesis which is always satisfied by the nerve fiber models considered
in \cite{ChBe-94}, \cite{GrSl-85}). In \cite{ZaZa-05} we have already proved the existence of
threshold values, associated to a general shape of $n(x)$ which imply, for the
particular case of $(\ref{enne}),$ that the only periodic solution is the trivial one
when $n_0$ and $\alpha$ are sufficiently small, while we proved the existence of at least
two nontrivial $\beta$-periodic solutions (under a weak technical assumption on $F(s)$)
when $n_1$ is sufficiently large and $\alpha$ is close to $\beta.$
Now, in view of Theorem \ref{th-5.1} and Theorem \ref{th-5.3} we propose the splitting
\[
n(x) = n_1 + {\tilde n}(x)
\]
where
\[
{\tilde n}(x) = \;
\left\{
\begin{array}{ll}
0\quad \mbox{if }\, x\in \,]0,\alpha[\, \mod \beta\\
n_0 - n_1\quad \mbox{if }\, x\in \,]\alpha,\beta[\, \mod \beta\\
\end{array}
\right.
\]
By this choice ${\bar n} = n_1$ and, moreover, we have
\[
|{\tilde n}|_1 = |{\tilde n}|_{L^1([0,m\beta])} = m(n_1 - n_0)(\beta - \alpha).
\]
In spite of the fact that we take ${\bar n} = n_1$ large, we are allowed to make
$|{\tilde n}|_1$ as small as we like, by taking $\alpha$ sufficiently close to $\beta.$
Hence Theorem \ref{th-5.1} and Theorem \ref{th-5.3} can be applied and the existence
of a large number of harmonic and ``true'' subharmonic solutions for the equation $(\ref{eq-1.1})$
is guaranteed.
}
\end{example}

\hrule

\end{document}